\documentclass[10pt]{article}
\usepackage{blindtext}
\usepackage[margin=1in]{geometry} 
\usepackage{amsmath,amsthm,amssymb, graphicx, multicol, array, mathtools}
\newtheorem{theorem}{Theorem}

\newtheorem{remark}{Remark}[theorem]
\newtheorem{example}{Example}[theorem]
\newtheorem{lemma}[theorem]{Lemma}
\usepackage[nottoc]{tocbibind}
\theoremstyle{definition}
\newtheorem{definition}{Definition}[section]
\usepackage{fancyhdr}
\usepackage{listings}

\usepackage{indentfirst}
\usepackage{color}
\usepackage[colorlinks=true, urlcolor=cyan]{hyperref}

\title{\textbf{A System of PDEs for the Baik-Rains Distribution}	}
\author{Xincheng Zhang\\xincheng.zhang@mail.utoronto.ca\\  Department of Mathematics\\
	University of Toronto}

\begin{document}
	\maketitle
	\begin{abstract}
		It has been discovered that the Kadomtsev-Petviashvili(KP) equation governs the distribution of the fluctuation of many random growth models, in particular, the Tracy-Widom distributions appear as special self-similar solutions of the KP equation. We prove that the anti-derivative of Baik-Rains distribution, which governs the fluctuation of the models in the KPZ universality class starting with stationary initial data, satisfies the KP equation. We start from a determinantal formula of the generating function of the KPZ equation, which satisfies the KP. Then we observe that the equation still holds in the large time limit.
	\end{abstract}
	\section{Introduction}

	The fluctuation for the KPZ universality class depends on the initial data. Let $h(x,t)$ be the solution of Kardar-Parisi-Zhang equation:
	\begin{equation}
	\partial_th(x,t) = \frac{\lambda}{2}(\partial_xh(x,t))^2+\nu \partial_x^2h(x,t)+\sqrt{D}\xi(x,t)
	\end{equation}
	\noindent $\mu >0,\lambda, D\neq 0$ are fixed parameters. $\xi(x,t)$ is the Gaussian space-time white noise
	\begin{equation}
	\mathbb{E}(\xi(x,t)\xi(y,s))= \delta(x-y)\delta(t-s)
	\end{equation}
	The equation is ill-posed because the quadratic non-linear term cannot make sense for a realization of a solution. A typical solution $h(x,t)$ looks like a Brownian motion in variable $x$. One way to make sense of the equation is through Hopf-Cole transformation.
	The Hopf-Cole solution of the KPZ equation is defined to be: $h(t,x) = -\log z(t,x)$, where $z(t,x)$ is the solution of the stochastic heat equation with multiplicative white noise:
	\begin{equation}
	\partial_tz(t,x) = \frac{1}{2}\partial_x^2z(t,x) + z(t,x)\xi(t,x)
	\end{equation}
	which is well-posed interpreted as It$\hat{\text{o}}$ integral.\\
	\\
	\indent All the universal fluctuation behavior can be observed on large space and time scales under the KPZ scaling:
	\begin{equation}
	h_\epsilon(t,x) = \epsilon^{1/2}h(\epsilon^{-3/2}t,\epsilon^{-1}x)
	\end{equation} 
	The fluctuation depends on the initial condition. If the initial condition  is $z(0,x) = \delta_0$,which corresponds to $h(0,x) = -\infty$ if $x \neq 0$ and $h(0,x) = 0$ if $x = 0$, i.e. KPZ starting from narrow wedge initial condition, we observe that when $t \rightarrow \infty$:
	\begin{equation}
	-2^{1/3}t^{-1/3}(h(t,2^{1/3}t^{2/3}x)-2^{-1/3}t^{1/3}x^2-\frac{t}{24}-\log\sqrt{2\pi t}) \rightarrow \mathcal{A}_2(x)
	\end{equation}
	If the initial condition is $z(0,x) = 1$, we observe:
	\begin{equation}
	-2^{1/3}t^{-1/3}(h(t,2^{1/3}t^{2/3}x)-\frac{t}{24}) \rightarrow \mathcal{A}_1(x)
	\end{equation}
	If the initial condition is $z(0,x) = e^{B(x)}$, where $B(x)$ is the two-sided Brownian motion with $B(0) = 0$, we observe:
	\begin{equation}
	-2^{1/3}t^{-1/3}(h(t,2^{1/3}t^{2/3}x)-\frac{t}{24}) \rightarrow \mathcal{A}_{\text{stat}}(x)
	\end{equation}
	$\mathcal{A}_1(x),\mathcal{A}_2(x),\mathcal{A}_{\text{stat}}(x)$ are stochastic processes whose finite dimensional distributions are given by Fredholm determinant. These are the conjectured processes which govern the long-time fluctuation of models which belong to the KPZ universality class. $\mathcal{A}_1(x)$ is a stationary process, whose one-point distribution is the Tracy-Widom GOE distribution. The one point marginal of $\mathcal{A}_2(x)-x^2$ is given by the Tracy-Widom GUE distribution. The one point marginal of $\mathcal{A}_{\text{stat}}(x)$ is given by the Baik-Rains distribution\\
	\\
	\indent In the paper \cite{kp}, it was stated that the GUE and GOE Tracy-Widom distributions are seen to arise as special similarity solutions of the scalar Kadomtsev-Petviashvili(KP) equation:
	\begin{equation}\label{kpequation}
	\partial_t\phi+\phi\partial_r\phi+\frac{1}{12}\partial_r^3\phi+\frac{1}{4}\partial_r^{-1}\partial_x^2\phi = 0
	\end{equation} In this paper, we explain how the Baik-Rains distribution can be seen as a similarity solution of the KP equation. We first explain how the GUE and GOE distributions arise as similarity solutions of the KP equation, then give the definition of Baik-Rains distribution and then state our main results.\\

	\begin{example}
		\textbf{Tracy-Widom GUE distribution}\cite{kp}:
		If we look for a self-similar solution of (\ref{kpequation}) of the form
		\begin{equation}
		\phi(t,x,r) = t^{-2/3}\psi(t^{-1/3}r+t^{-4/3}x^2)
		\end{equation}
		This turns (\ref{kpequation}) to 
		\begin{equation}
		\psi'''+12\psi\psi'-4r\psi'-2\psi = 0
		\end{equation}
		If we look for solutions of the form $\psi = -q^2$, then the above equation becomes the Painl\'eve II equation:
		\begin{equation}
		q'' = rq+2q^3
		\end{equation} in this way we just recover the GUE distribution, since $F_{\text{GUE}} = \exp\{-\int_s^{-\infty}du(u-s)q^2(u)\}$\\
	\end{example}
	\begin{example}
			\textbf{Tracy-Widom GOE distribution}\cite{kp}:
		If we look for a self-similar solution of (8) of the form
		\begin{equation}
		\phi(t,r) = (t/4)^{-2/3}\psi((t/4)^{-1/3}r)
		\end{equation}
		then the above equation becomes:
		\begin{equation}
		\psi'''+12\psi'\psi-r\psi'-2\psi=0
		\end{equation}
		If we look for solutions of form $\psi = \frac{1}{2}(q'-q^2)$, we get the Painl\'eve II equation again, thus we recover the 
		\begin{equation}
		F_{\text{GOE}}(r) = \exp\{-\frac{1}{2}\int_r^{\infty}q(u)du\}F_{\text{GUE}}(r)^{1/2}
		\end{equation}
	\end{example}

	Now let's look at the definition of the Baik-Rains distribution.
	\begin{definition}\cite{scaling}
		Ai is the Airy function. For $\tau, s \in \mathbb{R}$, define
		\begin{gather}
		\hat{\Phi}_{w,s}(x)=\int_{\mathbb{R}_{-}}dze^{wz}K_{\text{Ai,s}}(z,x)e^{ws}\\
		\hat{\Psi}_{w,s}(y) =\int_{\mathbb{R}_{-}}dze^{wz}\text{Ai}(y+z+s)\\
		\rho_s(x,y) = (\boldsymbol{I}-P_0\boldsymbol{K}_{\text{Ai,s}}P_0)^{-1}(x,y) 
		\end{gather}
		\noindent $P_0(x) = \boldsymbol{I}_{x\geq 0}$ is the projection operator, the shifted Airy kernel
		\begin{equation}
		\boldsymbol{\hat{K}_{\text{Ai,s}}}(x,y) = \int_0^\infty d\lambda\text{Ai}(x+\lambda+s)\text{Ai}(y+\lambda+s)
		\end{equation}
		the Tracy-Widom GUE distribution can be written as 
		\begin{equation}
		F_{\text{GUE}}(s)=\det(\boldsymbol{I}-\boldsymbol{P}_0\boldsymbol{K}_{\text{Ai,s}}\boldsymbol{P}_0)
		\end{equation}
		then we define the function $g(s,w)$ which appear as a component in the Baik-Rains distribution
		\begin{equation}
		g(s,w)  = e^{-\frac{1}{3}w^3}[\int_{\mathbb{R}^2_{-}}dxdye^{w(x+y)}\text{Ai}(x+y+s)+\int_{\mathbb{R}^2_{+}}dxdy\hat{\Phi}_{w,s}(x)\rho_s(x,y)\hat{\Psi}_{w,s}(y)]
		\end{equation}
		Finally the Baik-Rains distribution is defined to be
        \begin{equation}
			F_\tau(r) = \frac{\partial}{\partial r}(
			g(r+\tau^2,\tau/2)F_{\text{GUE}}(r+\tau^2))
			\end{equation}

	\end{definition}
 Notice that the function $g(s, w)$ is also derived in solving the PNG model, using the Riemann-Hilbert technique, given as solutions of a set of differential equations. We also present this equivalent definition here\cite{limiting}: Let $u(x)$ be the solution of the Painl\'eve II equation
	\begin{equation}
	u_{xx} = 2u^3+xu
	\end{equation}
	with the boundary condition
	\begin{equation}\label{painleve2}
	u(x) \sim -\text{Ai}(x) \text{ as } x \rightarrow+\infty
	\end{equation}
	$v(x)$ is defined to be
	\begin{equation}
	v(x)= \int_{\infty}^x(u(s))^2ds
	\end{equation}
	then the Tracy-Widom distribution can be defined in terms of $u$ and $v$. Set
	\begin{equation}
	F(x) = \exp(\frac{1}{2}\int_x^{\infty}v(x)ds) = \exp(-\frac{1}{2}\int_x^\infty(s-x)(u(s))^2ds),\qquad	E(x) = \exp(\frac{1}{2}\int_x^{\infty}u(s)ds)
	\end{equation}
	then
	\begin{gather}
	F_{\text{GUE}}(x) = F(x)^2 =  \exp(\int_x^\infty(s-x)(u(s))^2ds)\\
	F_{\text{GOE}}(x)= F(x)E(x)
	\end{gather}
	then $F_\tau(r)$ is defined to be $F_\tau(r) = H(r+\tau^2;\tau/2,-\tau/2)$, where
	\begin{equation}
    H(x;w,-w) = \{y'(x,w)-y(x,w)v(x))\}F_{GUE}(x) = \partial_x(y(x,w)F_{\text{GUE}}(x))
    \end{equation}
	where
	\begin{equation}\label{functiony}
	y(x,w) = (2u^2+x-4w^2)a(x;w)a(x;-w)-(u'+2wu)b(x;w)a(x;-w)-(u'-2wu)a(x;w)b(x;-w)
	\end{equation}
	Here functions $a(x;w),b(x;w)$ arise in the Painl\'eve II Riemann-Hilbert problems. In this paper, we do not need the precise definition of $a(x;w),b(x;w)$, thus we skip the definition here. What we need is the following identities:
    \begin{align}\label{identities}
		\partial_x a(x,w) &= u(x)b(x,w)\\
	\partial_x b(x,w) &= u(x)a(x,w) - 2wb(x,w)\\
	\partial_w a(x,w) &= 2u(x)^2a(x,w) - (4wu(x)+2u'(x))b(x,w)\\
	\partial_w b(x,w) &= (-4wu(x)+2u'(x))a(x,w) + (8w^2-2x-2u(x)^2)b(x,w)\\
	a(x,w) &= -b(x,-w)e^{\frac{8}{3}w^3-2wx}\\\label{identitieslast}
	b(x,w) &= -a(x,-w)e^{\frac{8}{3}w^3-2wx}
	\end{align}

	\noindent It is proven in \cite{scaling} that $y(s,w) = g(s, w)$, thus two definitions we presented are the same. In the following part of the paper, we will identify functions $g(s,w)$ and $y(s,w)$.\\
	\\
	\indent The main result we discovered about $F_\tau(r)$ is:
	\begin{theorem}\label{thm1}		
	    $F_\tau(r)$ is defined to be the partial derivative of $
		g(r+\tau^2,\tau/2)F_{\text{GUE}}(r+\tau^2)$ in $r$. If we consider certain scaling form of this anti-derivative of $F_\tau(r)$, which is\begin{equation}
			(\partial_r^{-1}F_{t^{-2/3}x})(t^{-1/3}r)= F_{\text{GUE}}(t^{-1/3}r+t^{-4/3}x^2)g(t^{-1/3}r+t^{-4/3}x^2,t^{-2/3}x/2)
		\end{equation} its logarithmic derivative $\phi_{br}(x,t,r) = \partial_r^2\log(\partial_r^{-1}F_{t^{-2/3}x})(t^{-1/3}r)$, satisfies the KP equation:
		\begin{gather}
		\partial_t\phi_{br} + \phi\partial_r\phi_{br}+\frac{1}{12}\partial_r^3\phi_{br}+\frac{1}{4}\partial^{-1}_r\partial_x^2
		\phi_{br} = 0
		\end{gather}
	\end{theorem}
    \begin{remark}
    	Now we know $(\partial_r^{-1 }F_{t^{-2/3}x})(t^{-1/3}r) = F_{\text{GUE}}(t^{-1/3}r+t^{-4/3}x^2)g(t^{-1/3}r+t^{-4/3}x^2,t^{-2/3}x/2)$ and $ \partial_r^2\log F_{\text{GUE}}$ also satisfies the KP equation, we denote it as $\phi_{gue} = \partial_r^2\log F_{\text{GUE}} $, so the Thm. 1 is equivalent to that $\psi(x,t,r) = \partial_r^2\log g$ satisfies the following equation:
    	\begin{gather} \label{modifiedkp}
    		\partial_t\psi + \psi\partial_r\psi+\frac{1}{12}\partial_r^3\psi+\frac{1}{4}\partial^{-1}_r\partial_x^2
    		\psi+\phi_{gue}\partial_r\psi+\psi\partial_r\phi_{gue} = 0
    	\end{gather}
    	In this paper, we will prove that $g$ satisfies (\ref{modifiedkp}), which implies Thm. 1.
    \end{remark}
    \noindent There are two ways we will show that Thm. \ref{thm1} is true. One way is to directly substitute equation (\ref{functiony}) into the equation (\ref{modifiedkp}). Using the identities in (\ref{identities})-(\ref{identitieslast}), we found that equation (\ref{modifiedkp}) holds. The other way is that by observing the Baik-Rains distribution govern the fluctuation of large time scale of the KPZ equation, we use the results that the generating function of the KPZ equation satisfies KP \cite{doussal}, taking limit in $t$, obtaining Thm. \ref{thm1}. But part of the second method is still non-rigorous.\\
	\\

	\section{Proof of the Theorem 1}
	
	 A key fact we will use is the following theorem\cite{kp}:
	\begin{theorem}
		If a function can be written in the Fredholm determinant form, i.e. $F(x,t,r)= \det(\boldsymbol{I}-\boldsymbol{K})_{L^2[0,\infty)}$, and if the integral kernel $\boldsymbol{K}(u,v,x,t,r)$ satisfies the following relations:
		\begin{equation}
		\begin{aligned}\label{kernelrelation}
		\partial_r\boldsymbol{K} = (\partial_u+\partial_v)\boldsymbol{K}\\
		\partial_t\boldsymbol{K} = -\frac{1}{3}(\partial_u^3+\partial_v^3)\boldsymbol{K}\\
		\partial_x\boldsymbol{K} = (\partial_u^2-\partial_v^2)\boldsymbol{K}
		\end{aligned}
		\end{equation}
		then $\phi(x,t,r) = \partial_r^2\log F$ satisfy the scalar KP equation:
	\end{theorem}
	\begin{equation}
	\partial_t\phi + \phi\partial_r\phi+\frac{1}{12}\partial_r^3\phi+\frac{1}{4}\partial^{-1}_r\partial_x^2
	\phi = 0
	\end{equation}
	\begin{remark}
		The fact that equation (\ref{kernelrelation}) leads to the KP equation was discovered several times \cite{zs74}\cite{po89}, but its appearance in the context of random fluctuation interfaces was discovered in \cite{kp}
	\end{remark}
	\noindent It is checked in \cite{kp} \cite{doussal} that if $h(x,t)$ is the solution of the KPZ equation with the half-Brownian initial data or narrow wedge initial data, their generating functions $G(t,x,r) = \mathbb{E}[\exp\{-e^{h(t,x)+\frac{t}{12}-r}\}]$ can be written in the Fredholm determinant form, having a kernel which satisfying equation (\ref{kernelrelation}). It is also checked in \cite{doussal} that if $h(x,t)$ is the solution of the KPZ equation with initial condition to be the drifted Brownian motion, its generating function also satisfies the KP equation, by studying the moments of $e^{h(x,t)}$. Here we checked it using a different method, which agreed with the results in \cite{doussal}  \\
	\\
	\begin{theorem}\cite{height}
		Let $z_{b,\beta}(t,x)$ denote the solution to the stochastic heat equation with initial data $z(0,x) = \exp(B_{b,\beta}(x))$, where $B_{b,\beta}(x)$ is a two-sided Brownian motion with drift $\beta$ to the left of 0 and drift $b$ to the right of 0, with $\beta > b$, that is, $B_{b,\beta} = \boldsymbol{1}_{x\leq 0}(B^l(x)+\beta x)+\boldsymbol{1}_{x >0}(B^r(x)+bx)$, where $B^l:(-\infty,0] \rightarrow \mathbb{R}$ is a Brownian motion without drift pinned at $B^l(0)=0$ and $B^r:[0,\infty) \rightarrow \mathbb{R}$ is an independent Brownian motion pinned at $B^r(0) = 0$, then for $S >0$\\
		\begin{equation}\label{expectationK}
		\mathbb{E}[2(Se^{\frac{x^2}{2t}+\frac{t}{24}}z_{b,\beta}(t,x))^{\frac{\beta-b}{2}}K_{-(\beta-b)}(2\sqrt{Se^{\frac{x^2}{2t}+\frac{t}{24}}z_{b,\beta}(t,x)})] = \Gamma(\beta - b)\det(\boldsymbol{I}-\boldsymbol{K}_{b+\frac{x}{t},\beta+\frac{x}{t}})_{L^2(\mathbb{R}_{+})}
		\end{equation}
		where $K_\nu(z)$ is the modified Bessel function of order $\nu$ and the kernel on the right-hand side is given by\\
		\begin{equation}	
		\boldsymbol{K}_{b,\beta}(x,y) = \frac{1}{(2\pi i)^2}\int dw\int dz\frac{\sigma \pi S^{\sigma(z-w)}}{\sin(\sigma\pi(z-w))}\frac{e^{z^3/3-zy}}{e^{w^3/3-wx}}\frac{\Gamma(\beta-\sigma z)}{\Gamma(\sigma z-b)}\frac{\Gamma(\sigma w-b)}{\Gamma(\beta-\sigma w)}
		\end{equation}
		where\\
		\begin{equation}
		\sigma = (2/t)^{1/3}
		\end{equation}
		The integration contour for $w$ is from $-\frac{1}{4\sigma}-i\infty$ to $-\frac{1}{4\sigma}+i\infty$ and crosses the real axis between $b$ and $\beta$. The other contour for $z$ goes from $\frac{1}{4\sigma}-i\infty$ to $\frac{1}{4\sigma}+i\infty$, it also crosses the real axis between $b$ and $\beta$ and it does not intersect the contour for $w$
	\end{theorem}
	\noindent The left hand-side of (\ref{expectationK}) is exactly the generating function, since for $\alpha > a$, we have\cite{height}:\\
	\begin{equation}
	\mathbb{E}[2(uz(\tau,N))^{\frac{\alpha - a}{2}}K_{-(\alpha - a)}(2\sqrt{uz(\tau,N)})] = \Gamma(\alpha - a)\mathbb{E}[e^{-uz(\tau,N)\gamma}]
	\end{equation}
	where $\frac{1}{\gamma}$ has gamma distribution with parameter $\alpha - a$.\\
	Thus we have:
	\begin{equation}
	\mathbb{E}[\exp\{-Se^{\frac{x^2}{2t}+\frac{t}{24}+h_{b,\beta}(t,x)+\log \gamma}\}] = \det(\boldsymbol{I}-\boldsymbol{K}_{b+\frac{x}{t},\beta+\frac{x}{t}})_{L^2(\mathbb{R}_{+})}
	\end{equation}
	For $S = e^{\tau^2+r}$, where $\tau$ is related to $x,b,t$ as $x = bt+\frac{2\tau}{\sigma^2}$, and when $b = 0$, we observe that $\boldsymbol{K}_{b+\frac{x}{t},\beta+\frac{x}{t}}$  satisfy relation (\ref{kernelrelation}). The reason that we make this specific choice of $S$ and $\tau$ will be clear from the later context, this is the scaling that gives the Baik-Rains distribution. We will do some transformations on the integral kernel so that equation (\ref{kernelrelation}) becomes obvious while the operator remains unchanged. When $S = e^{\tau^2+r}$,
	\begin{align}
	\boldsymbol{K}_{b+\frac{x}{t},\beta +\frac{x}{t}}(u,v) = \frac{1}{(2\pi i)^2}\int dw\int dz\frac{\sigma \pi e^{-(z-w)(\tau^2+r)}}{\sin(\sigma \pi(z-w))}\frac{e^{z^3/3-zv}}{e^{w^3/3-wu}}\frac{\Gamma(\beta +\frac{x}{t}-\sigma z)}{\Gamma(\sigma z-b-\frac{x}{t})}\frac{\Gamma(\sigma w-b-\frac{x}{t})}{\Gamma(\beta+\frac{x}{t}-\sigma w)}
	\end{align}
	Now let $z = \sigma z - \frac{x}{t},w = \sigma w- \frac{x}{t} $, then 
	\begin{align}
	\boldsymbol{K}_{b+\frac{x}{t},\beta +\frac{x}{t}}(u,v) = \frac{1}{(2\pi i)^2}\int dw\int dz\frac{ \pi e^{-(z-w)(\tau^2+r)/\sigma}}{\sin( \pi(z-w))}\frac{e^{t(z^3+3z^2x/t+3zx^2/t^2)/6-\frac{1}{\sigma}(z+\frac{x}{t})v}}{e^{t(w^3+3w^2x/t+3wx^2/t^2)/6-\frac{1}{\sigma}(w+\frac{x}{t})u}}\frac{\Gamma(\beta - z)}{\Gamma( z-b)}\frac{\Gamma(w-b)}{\Gamma(\beta- w)}
	\end{align}
	Let $u = \frac{1}{\sigma}u, v= \frac{1}{\sigma}v$ and conjugate $e^{(v-u)x/t}$, and then change $r = \frac{1}{\sigma}r$ it becomes,
	\begin{align}
	\boldsymbol{K}_{b+\frac{x}{t},\beta +\frac{x}{t}}(u,v) = \frac{1}{(2\pi i)^2}\int dw\int dz\frac{ \pi e^{-(z-w)(b^2t/2+bx)}}{\sin( \pi(z-w))}\frac{e^{(tz^3+3z^2x)/6-z(v+r)}}{e^{(tw^3+3w^2x)/6-w(u+r)}}\frac{\Gamma(\beta - z)}{\Gamma( z-b)}\frac{\Gamma(w-b)}{\Gamma(\beta- w)}
	\end{align}
	When $b= 0$, the only term contain $x,t,r,u,v$ is $\frac{e^{(tz^3+3z^2x)/6-z(v+r)}}{e^{(tw^3+3w^2x)/6-w(u+r)}}$, which clearly satisfies equation (\ref{kernelrelation}). When $b \neq 0$, it will have extra terms  when differentiating $t,x$ coming from $e^{-(z-w)(b^2t/2+bx)}$,which fail the equation (\ref{kernelrelation}). The reason that we can directly take derivative under the integral sign is the following lemma:\\
	\begin{lemma}\cite{height}
		Let $f(z,\zeta)$ be a complex function in two variables and suppose that
		\begin{enumerate}
			\item $f$ is defined on $(z,\zeta)\in A \times C$ where $A$ is an open set and $C$ is a contour
			\item For each $z\in A$, define the contour $\gamma = \{z+re^{it}: 0\leq t\leq 2\pi\}$ with a sufficiently small r such that also the disc around z with radius r lies in A. Suppose that for each $z \in A$\\
			\begin{equation}
			\int_C\int_\gamma|f(u,\zeta)||du||d\zeta| < \infty
			\end{equation}
			\item For each $\zeta \in C, z \rightarrow f(z,\zeta)$ is analytic in A
			\item For each $z \in A, \zeta \rightarrow f(z,\zeta)$ is continuous on C
		\end{enumerate}
		Then \begin{equation}
		F(z) = \int_Cf(z,\zeta)d\zeta
		\end{equation}
		is analytic in A with $F'(z) = \int_C\frac{\partial}{\partial z}f(z,\zeta)d\zeta$
	\end{lemma}
	\noindent It can be easily seen that condition 2 is satisfied. Since $e^{z^3/3}$ decay along $C_z$ as $e^{-c|\text{Im}(z)|^2}$, $e^{w^3/3}$ decay along $C_w$ as $e^{-c|\text{Im}(w)|^2}$. Using gamma ratio formula, as $|z| \rightarrow \infty$,
	\begin{equation}
	\left| \frac{\Gamma(\beta -\sigma z)}{\Gamma(\sigma z-b)}\right| \simeq |z|^{\beta + b-2\sigma\text{Re}(z)}
	\end{equation}
	Similarly we have the same bound for  large $w$. Thus if we integrate $\beta$ on some finite contour, we have the same polynomial bounds. Thus the whole integrands decay exponentially on the contour, so condition (2) is satisfied.\\
	\\
	\indent In order to get the formula for stationary initial data, we need to take the limit as $\beta \rightarrow b$ and set $b= 0$. To do so, we need to rewrite the kernel $\boldsymbol{K}_{b,\beta}$. Two contours in the integral kernel of $\boldsymbol{K}_{b,\beta}$ intersect the real axis between the pole at $b/ \sigma$ and $\beta / \sigma$, so when $\beta \rightarrow b$, two contours will collide. Hence by using the residue theorem, we cross the pole at $b/ \sigma$ with the $w$ integration contour
	and cross the pole at $\beta / \sigma$ with the $z$ integration contour, both manipulations resulting in a residue term. We have\cite{height}

	\begin{align}
	\boldsymbol{K}_{b,\beta} &= \boldsymbol{\bar{K}}_{b,\beta}+q_{b,\beta}(x)r_{\beta}(y)\frac{1}{\sigma\Gamma(\beta-b)}+r_{-b}(x)q_{-\beta,-b}(y)\frac{1}{\sigma\Gamma(\beta-b)}\\
	&+\frac{\sigma\pi S^{\beta - b}}{\sin(\pi(\beta - b))}r_{-b}(x)r_{\beta}(y)\frac{1}{\sigma^2\Gamma(\beta - b)^2}
	\end{align}

	where
	\begin{align}
	\boldsymbol{\bar{K}}_{b,\beta} &= \frac{1}{(2\pi i)^2}\int_{-\frac{1}{4\sigma}+i\mathbb{R}} dw\int_{\frac{1}{4\sigma}+i\mathbb{R}} dz\frac{\sigma \pi S^{\sigma(z-w)}}{\sin(\sigma\pi(z-w))}\frac{e^{z^3/3-zy}}{e^{w^3/3-wx}}\frac{\Gamma(\beta-\sigma z)}{\Gamma(\sigma z-b)}\frac{\Gamma(\sigma w-b)}{\Gamma(\beta-\sigma w)}\\
	q_{b,\beta}(x)&= \frac{1}{2\pi i}\int_{-\frac{1}{4\sigma}+i\mathbb{R}}dw\frac{\sigma \pi S^{\beta-\sigma w}}{\sin(\pi(\beta-\sigma w))}e^{-w^3/3+wx}\frac{\Gamma(\sigma w-b)}{\Gamma(\beta- \sigma w)}\\
	r_b(x) &= e^{b^3/(3\sigma^3)-bx/\sigma}
	\end{align}
	\noindent Notice that the only difference between $\boldsymbol{K}_{b,\beta}$ and $\boldsymbol{\bar{K}_{b,\beta}}$ is that they have different contour. We can write
	\begin{equation}
	\boldsymbol{K_{b,\beta}} = \boldsymbol{\bar{K}_{b,\beta}}+\sum_{i=1}^3f_i(x)g_i(y)
	\end{equation}
	with suitable $f_i,g_i$, then for Fredholm determinant, we have the following formula
	\begin{equation}\label{midkernelidentity}
	\det(\boldsymbol{I}-\boldsymbol{K_{b,\beta}}) = \det(\boldsymbol{I}-\boldsymbol{\bar{K}_{b,\beta}})\det[\delta_{i,j}-\langle(\boldsymbol{I}-\boldsymbol{\bar{K}_{b,\beta}})^{-1}f_i,g_j\rangle]_{i,j=1}^3
	\end{equation}
	Here $\boldsymbol{\bar{K}}_{b,\beta}$ also satisfies relation (\ref{kernelrelation}) if $S = e^{\tau^2+r}$, because the only difference between $\boldsymbol{\bar{K}}_{b,\beta}$ and $\boldsymbol{K}_{b,\beta}$ is that they have different contour, for which the lemma 4 still apply. We denote  $\phi = \partial_r^2\log\det(\boldsymbol{I}-\boldsymbol{K_{b,\beta}}),\bar{\phi} = \partial_r^2\log\det(\boldsymbol{I}-\boldsymbol{\bar{K}_{b,\beta}})$ and $\alpha = \partial_r^2\log \det[\delta_{i,j}-\langle(\boldsymbol{I}-\boldsymbol{\bar{K}_{b,\beta}})^{-1}f_i,g_j\rangle]_{i,j=1}^3$, we have $\phi = \bar{\phi} + \alpha$, and both $\phi, \bar{\phi}$ satisfy the KP equation\\
	\begin{gather}
	\partial_t\bar{\phi} + \bar{\phi}\partial_r\bar{\phi} + \frac{1}{12}\partial_r^3\bar{\phi}+\frac{1}{4}\partial_r^{-1}\partial_x^2\bar{\phi}=0\\
	\partial_t(\bar{\phi}+\alpha) + (\bar{\phi}+\alpha)\partial_r(\bar{\phi}+\alpha) + \frac{1}{12}\partial_r^3(\bar{\phi}+\alpha)+\frac{1}{4}\partial_r^{-1}\partial_x^2(\bar{\phi}+\alpha)=0
	\end{gather}
	combine these two equations, we obtain an equation for $\alpha$ which involves $\bar{\phi}$:\\
	\begin{equation}
	\partial_t\alpha + \alpha\partial_r\alpha + \frac{1}{12}\partial_r^3\alpha+\frac{1}{4}\partial_r^{-1}\partial_x^2\alpha+\alpha\partial_r\bar{\phi}+\bar{\phi}\partial_r\alpha=0
	\end{equation}
	This is the first place where we obtain equation (\ref{modifiedkp}). Now we want to see how equation (\ref{midkernelidentity}) leads us to the Baik-Rains distribution.\\
	\\
	\\
	In the limit as $\beta \rightarrow b$, we have the following results
	\begin{theorem}\cite{height}
		Let $b+\frac{x}{t} \in (-\frac{1}{4},\frac{1}{4})$ be fixed. For the kernel $\boldsymbol{K}_{b,\beta}$, we have
		\begin{equation}\label{xilimit}
		\lim_{\beta \rightarrow b}\frac{1}{\beta - b}\det(\boldsymbol{I}-\boldsymbol{K}_{b+\frac{x}{t},\beta+\frac{x}{t}}) = \frac{1}{\sigma}\Xi(S,b+\frac{x}{t},\sigma)
		\end{equation}
		where\\
		\begin{equation}
		\begin{aligned}
		\Xi(S,b,\sigma) &= -\det(\boldsymbol{I}-\boldsymbol{\bar{K}_b})[\frac{b^2}{\sigma^2}+\sigma(2\gamma_E+\ln S)\\&+
		\langle(\boldsymbol{I}-\boldsymbol{\bar{K}_b})^{-1}(\boldsymbol{\bar{K}_b}r_{-b}+q_b),r_b\rangle+\langle(\boldsymbol{I}-\boldsymbol{\bar{K}_b})^{-1}(r_{-b}+q_b),q_{-b}\rangle]
		\end{aligned}
		\end{equation}
		Where $\gamma_E$ represents Euler-Mascheroni constant, and for $b+\frac{x}{t} \in (-\frac{1}{4},\frac{1}{4}),\boldsymbol{\bar{K}}_{b+\frac{x}{t}} = \boldsymbol{\bar{K}}_{b+\frac{x}{t},b+\frac{x}{t}}, q_{b+\frac{x}{t}} = q_{b+\frac{x}{t},b+\frac{x}{t}}$\\
	\end{theorem}
	\noindent For simplicity, we will define
	\begin{gather}
	A(S,b,\sigma) = \det(\boldsymbol{I}-\boldsymbol{\bar{K}_b})\\
	B(S,b,\sigma) =[\frac{b^2}{\sigma^2}+\sigma(2\gamma_E+\ln S)+
	\langle(\boldsymbol{I}-\boldsymbol{\bar{K}_b})^{-1}(\boldsymbol{\bar{K}_b}r_{-b}+q_b),r_b\rangle+\langle(\boldsymbol{I}-\boldsymbol{\bar{K}_b})^{-1}(r_{-b}+q_b),q_{-b}\rangle]\\\label{phi}
	\phi(x,t,r) = \partial_r^2\log A(S,b+\frac{x}{t},\sigma)\\\label{psi}
	\psi(x,t,r) = \partial^2_r \log B(S,b+\frac{x}{t},\sigma)
	\end{gather} .\\
	Since formula (\ref{xilimit}) is obtained by taking limits of right-hand side of equation (\ref{midkernelidentity}), thus we speculate that the same system of PDEs will also hold for (\ref{xilimit}). What we need to show is that the limit of a sequence of functions satisfy the same PDE.\\
	We have the following results about the Fredholm determinant:
	\begin{lemma}\cite{height}
		The Fredholm determinant $\det(\boldsymbol{I}-\boldsymbol{K_{b,\beta}}), \det(\boldsymbol{I}-\boldsymbol{\bar{K}_{b,\beta}})$ are analytic functions of the parameter $b$ and $\beta$ as long as $b<\beta$
	\end{lemma}
	\noindent The reason that we require $b < \beta$ simply because otherwise the kernel $\boldsymbol{K}_{b,\beta}$ is not well defined.\\
	we want to show in the limit, two parts of right-hand side of (\ref{xilimit})  are both analytic. For the first part, we show that $\det(\boldsymbol{I}-\boldsymbol{\bar{K}_{b+\frac{x}{t},\beta+\frac{x}{t}}}) \rightarrow \det(\boldsymbol{I}-\boldsymbol{\bar{K}_{b+\frac{x}{t},b+\frac{x}{t}}})$ uniformly on compact set. Since
	\begin{equation}
	\det(\boldsymbol{I}-\boldsymbol{\bar{K}_{b+\frac{x}{t},\beta+\frac{x}{t}}}) = \sum_{n=0}^{\infty}\frac{(-1)^n}{n!}\int_0^\infty\cdots\int_0^\infty dx_1\cdots dx_n\det[\bar{\boldsymbol{K}}_{b,\beta}(x_i,x_j)]_{i,j=1}^n
	\end{equation}
	We need to show $\bar{\boldsymbol{K}}_{b,\beta} \rightarrow \bar{\boldsymbol{K}}_{b,b}$ uniformly on compact set. That simply because for the integrand of the kernel, $\frac{\Gamma(\beta - \sigma z)}{\Gamma(\beta-\sigma w)} \rightarrow \frac{\Gamma(b - \sigma z)}{\Gamma(b-\sigma w)}$ uniformly on compact set. For the other part, $[\frac{b^2}{\sigma^2}+\sigma(2\gamma_E+\ln S)+
	\langle(\boldsymbol{I}-\boldsymbol{\bar{K}_b})^{-1}(\boldsymbol{\bar{K}_b}r_{-b}+q_b),r_b\rangle+\langle(\boldsymbol{I}-\boldsymbol{\bar{K}_b})^{-1}(r_{-b}+q_b),q_{-b}\rangle]$, the only part that is non-trivial is whether $(\boldsymbol{I}-\bar{\boldsymbol{K}}_{b})^{-1}$ is analytic. We need the following theorem which gives an explicit formula for the inverse which can be found\cite{peter}:
	\begin{theorem}
		Let $\boldsymbol{K}$ be a continuous kernel, let $D = \det(\boldsymbol{I}+\boldsymbol{K})$, and suppose that $D \neq 0$. Then the operator $\boldsymbol{I}+\boldsymbol{K}$ is invertible, and its inverse is $\boldsymbol{I}-D^{-1}\boldsymbol{R}$, where the integral kernel of $\boldsymbol{R}$ is
		\begin{equation}
		R(x,y) = \sum_0^{\infty}\frac{1}{k!}\int\cdots\int K\begin{pmatrix}
		x, & x_1, & \cdots, & x_k \\
		y, & x_1, & \cdots, & x_k
		\end{pmatrix}dx_1\cdots dx_k
		\end{equation} 
	\end{theorem}
	\noindent From Thm. 7 we can conclude that if $K(x,y)$ is an analytic function, so is $(\boldsymbol{I}-\boldsymbol{K})^{-1}$.  
	\noindent Thus we conclude without fully justification that $\phi(x,t,r) = \partial_r^2\log A(S,b+\frac{x}{t},\sigma)$ satisfies equation (\ref{kpequation}) and $\psi(x,t,r) = \partial^2_r \log B(S,b+\frac{x}{t},\sigma)$ satisfies equation (\ref{modifiedkp}). In order to make the argument complete, we need to prove that the limit of partial derivatives converge uniformly.
	\\
	
	\indent In the large time limit, we have the following results\cite{height}:
	\begin{equation}
	\lim_{t \rightarrow \infty}\Xi(e^{-\frac{\tau^2+r}{\sigma}},\tau\sigma,\sigma) = 	g(r+\tau^2,\tau/2)F_{\text{GUE}}(r+\tau^2)
	\end{equation}
	Here the relation of $\tau$ to $x,t$ is:
	\begin{equation}
	x=  -bt +\frac{2\tau}{\sigma^2}, \sigma  =(\frac{2}{t})^{1/3}
	\end{equation}
	The reason for having this scaling is that under this scaling, we have the following result:
	\begin{theorem}\cite{height}
		Let $b\in (-\frac{1}{4},\frac{1}{4})$ be fixed and consider any $\tau \in \mathbb{R}$. Define $\sigma = (2/t)^{1/3}$ and consider the scaling $x=  -bt +\frac{2\tau}{\sigma^2}$. Then, for any $r \in \mathbb{R}$
		
		\begin{gather}
		\lim_{t\rightarrow \infty}\mathbb{P}(\frac{h_b(t,x)+\frac{t}{24}(1+12b^2)-2^{1/3}b\tau t^{2/3}}{(t/2)^{1/3}} \leq r) = F_\tau(r)
		\end{gather}
	\end{theorem}
	\noindent When $b=0$, we have $\sigma^2x = 2\tau$. We can see that as $t \rightarrow \infty$, $x$ must goes to the infinity in the speed $x \sim t^{2/3} $ so that $\tau$ can be a meaningful number. For this reason, we want to take the limit in the following way: we rewrite all the results for $h_\epsilon(t,x) = \epsilon^{1/2}h(\epsilon^{-1}x,\epsilon^{-3/2}t)$, the equation scales to:
	\begin{equation}
	\partial_t h_\epsilon = \frac{1}{2}(\partial_x h)^2+\epsilon^{1/2}\frac{1}{2} \partial_x^2 h+\epsilon^{1/4}\xi
	\end{equation}
	Then we plug $\epsilon^{-1}x \rightarrow x,\epsilon^{-3/2}t \rightarrow t, \epsilon^{-1/2}r \rightarrow r$ into $\Xi(e^{-\frac{\tau^2+r}{\sigma}},\tau\sigma,\sigma)$, which we denoted as $\Xi_\epsilon$ then take $\epsilon$ goes to 0, we have 
	\begin{equation}
	\lim_{\epsilon \rightarrow 0}\Xi_\epsilon(e^{-\frac{\tau^2+r}{\sigma}},\tau\sigma,\sigma) = g(\sigma r+\tau^2,\tau/2)F_{\text{GUE}}(\sigma r+\tau^2)
	\end{equation}
	Now $\Xi$ is a function involve $\epsilon$. We write $\phi^\epsilon(x,t,r) = \phi(\epsilon^{-1}x,,\epsilon^{-\frac{3}{2}}t,\epsilon^{-1/2}r),\psi^\epsilon(x,t,r) = \psi(\epsilon^{-1}x,,\epsilon^{-\frac{3}{2}}t,\epsilon^{-1/2}r)$, where $\phi,\psi$ are defined in (\ref{phi}),(\ref{psi}). Equation (\ref{kpequation}), (\ref{modifiedkp}) re-scale as follows:
	\begin{gather}
	\epsilon^{5/2}\partial_t\phi^\epsilon +\epsilon^{5/2}\phi^\epsilon\partial_r\phi^\epsilon + \frac{1}{12}\epsilon^{5/2}\partial_r^3\phi^\epsilon+\epsilon^{5/2}\frac{1}{4}\partial_r^{-1}\partial_x^2\phi^\epsilon=0\\
	\epsilon^{5/2}\partial_t\psi^\epsilon + \epsilon^{5/2}\psi^\epsilon\partial_r\psi^\epsilon + \epsilon^{5/2}\frac{1}{12}\partial_r^3\psi^\epsilon+\epsilon^{5/2}\frac{1}{4}\partial_r^{-1}\partial_x^2\psi^\epsilon+\epsilon^{5/2}\psi^\epsilon\partial_r\phi^\epsilon+\epsilon^{5/2}\phi^\epsilon\partial_r\psi^\epsilon=0
	\end{gather}
	Thus $\phi^\epsilon,\psi^\epsilon$ satisfy equation (\ref{kpequation}),(\ref{modifiedkp}). We conjecture without fully justification that in the limit as $\epsilon \rightarrow 0$, $\phi^\epsilon \rightarrow \partial^2_r\log F_{\text{GUE}}(\sigma r+\tau^2)$ satisfies (\ref{kpequation}), $\psi^\epsilon \rightarrow \partial^2_r\log g(\sigma r+\tau^2,\tau/2)$ satisfies (\ref{modifiedkp}), which is our Thm. \ref{thm1}. \\
	\\
	
	\section{Direct Verification}
	
	Now we are going to directly check that Let $B =  \partial^2_r\log y( t^{-1/3}r+t^{-4/3}x^2,\frac{1}{2}t^{-2/3}x)$, we use $y'$ to denote $\partial_1 y$ which is the partial derivative with respect to the first variable. $\phi = \partial_r^2\log F_{\text{GUE}}(t^{-1/3}r+t^{-4/3}x^2) = -t^{-2/3}u^2$, $u$ is defined in (\ref{painleve2}).
	\begin{align}
	\partial_rB &= t^{-3/3}(\tfrac{y'''}{y}-\tfrac{3y'y''}{y^2}+\tfrac{2y'^3}{y^3})\\
	B\partial_rB &= t^{-5/3}(\tfrac{y'''y''}{y^2}-\tfrac{3y'y''^2}{y^3}+\tfrac{5y'^3y''}{y^4}-\tfrac{y'''y'^2}{y^3}-\tfrac{2y'^5}{y^5})\\
	\phi\partial_rB &= t^{-5/3}(-u^2)(\tfrac{y'''}{y}-\tfrac{3y'y''}{y^2}+\tfrac{2y'^3}{y^3})\\
	B\partial_r\phi &= t^{-5/3}(-2uu')(\tfrac{y''}{y}-\tfrac{y'^2}{y^2})\\
	\partial_r^2B &=t^{-4/3}(\tfrac{y^{(4)}}{y} -\tfrac{3y''y''}{y^2}-\tfrac{4y'y'''}{y^2}+\tfrac{12y'^2y''}{y^3}-\tfrac{6y'^4}{y^4})\\
	\partial_r^3B &=t^{-5/3}(\tfrac{y^{(5)}}{y}-\tfrac{10y''y'''}{y^2}-\tfrac{5y'y^{(4)}}{y^2}+\tfrac{20y'^2y'''}{y^3}+\tfrac{30y'y''^2}{y^3}-\tfrac{60y'^3y''}{y^4}+\tfrac{24y'^5}{y^5})\\
	\partial_tB &= -\tfrac{2}{3}t^{-5/3}(\tfrac{y''}{y}-(\tfrac{y'}{y})^2)-t^{-2/3}(\tfrac{y'''}{y}-\tfrac{3y'y''}{y^2}+\tfrac{2y'^3}{y^3})(\tfrac{1}{3}t^{-4/3}r+\tfrac{4}{3}t^{-7/3}x^2)\\
	&- t^{-2/3}(\tfrac{\partial_2y''}{y}-\tfrac{y''\partial_2y}{y^2}-\tfrac{2y'\partial_2y'}{y^2}+\tfrac{2y'^2\partial_2y}{y^3})\tfrac{1}{3}t^{-5/3}x\\
	\partial_r^{-1}\partial_x^2B  &= 2t^{-5/3}(\tfrac{y''}{y}-\tfrac{y'^2}{y^2})+4x^2t^{-9/3}(\tfrac{y'''}{y}-\tfrac{3y'y''}{y^2}+\tfrac{2y'^3}{y^3})+\tfrac{1}{4}t^{-5/3}\partial_2^2(\tfrac{y'}{y})\\
	+&2xt^{-7/3}(\tfrac{\partial_2y''}{y}-\tfrac{y''\partial_2y}{y^2}-\tfrac{2y'\partial_2y'}{y^2}+\tfrac{2y'^2\partial_2y}{y^3})
	\end{align}
	Now plugin every term into equation (\ref{modifiedkp}), it becomes:
	\begin{align}
	&\partial_tB +B\partial_rB + \tfrac{1}{12}\partial_r^3B+\tfrac{1}{4}\partial_r^{-1}\partial_x^2B + \bar{\phi}\partial_r B+ B\partial_r \bar{\phi} \\ &= t^{-5/3}(\tfrac{1}{12}(\tfrac{y^{(5)}}{y}-\tfrac{5y'y^{(4)}}{y^2}+\tfrac{2y''y'''}{y^2}-\tfrac{6y'y''^2}{y^3}+\tfrac{8y'^2y'''}{y^3})-\tfrac{1}{6}(\tfrac{y''}{y}-\tfrac{y'^2}{y^2})\\
	&+\tfrac{1}{16}(\tfrac{\partial^2_2y'}{y}-\tfrac{2\partial_2y'\partial_2y}{y^2}-\tfrac{y'\partial_2^2y}{y^2}+\tfrac{2y'(\partial_2y)^2}{y^3})+\tfrac{1}{3}w(\tfrac{\partial_2y''}{y}-\tfrac{y''\partial_2y}{y^2}-\tfrac{2y'\partial_2y'}{y^2}+\tfrac{2y'^2\partial_2y}{y^3})\\
	&-\tfrac{1}{3}x(\tfrac{y'''}{y}-\tfrac{3y'y''}{y^2}+2\tfrac{y'^3}{y^3})-u^2(\tfrac{y'''}{y}-\tfrac{3y'y''}{y^2}+2\tfrac{y'^3}{y^3})-2uu'(\tfrac{y''}{y}-\tfrac{y'^2}{y^2}))
	\end{align}
	multiplied by $y^3t^{5/3}$,right-hand side becomes
	\begin{align}
	\label{equationsofy}
    &\tfrac{1}{12}(y^{(5)}y^2-5y'y^{(4)}y+2y''y'''y-6y'y''^2+8y'^2y''')-\tfrac{1}{6}(y''y^2-y'^2y)\\
	&+\tfrac{1}{16}(\partial^2_2y'y^2-2\partial_2y'\partial_2yy-y'\partial_2^2yy+2y'(\partial_2y)^2)+\tfrac{1}{3}w(\partial_2y''y^2-y''\partial_2yy-2y'\partial_2y'y+2y'^2\partial_2y)\\
	&-\tfrac{1}{3}x(y'''y^2-3y'y''y+2y'^3)-u^2(y'''y^2-3y'y''y+2y'^3)-2uu'(y''y^2-y'^2y)
	\end{align}
	Then we compute the derivative of $y$, in the following expressions, if we omit the variables, then it just means that's the variable in the definition; $a(-w), b(-w)$ represents $a(x,-w),b(x,-w)$:
	\begin{align}
	y &= (2u^2+x-4w^2)aa(-w)-(u'+2wu)ba(-w)-(u'-2wu)ab(-w)\\
	y' &= aa(-w)\\
	y'' &= uba(x,-w)+uab(-w)\\
	y''' &= (u'-2wu)ba(-w)+4u^2aa(-w)+(u'+2wu)ab(-w)\\
	y^{(4)} &=12uu'aa(-w)+(4u^3+u''+4wu'+4w^2u)ab(-w)\\
	&+(u''+4u^3-4wu'+4w^2u)ba(-w)\\
	y^{(5)} &= (12u'^2+16uu''+16u^4+16w^2u^2)aa(-w)\\
	&+ (24u^2u'+u'''+6wu''+12w^2u'+8wu^3+8w^3u)ab(-w)\\
	&+(24u^2u'+u'''-6wu''-8wu^3+12w^2u'-8w^3u)ba(-w)\\
	\partial_2y &= -8waa(-w)+2uab(-w)+(-2u)ba(-w)\\
	\partial_2^2y  &= (-8-16uu')aa(-w)+(16w^2u-16wu'+8u^3+4ux)ab(-w)\\
	&+ (16w^2u+16wu'+8u^3+4ux)ba(-w)
	\end{align}
	Observe that the derivative of $a(x,w),b(x,w)$ behave similarly like $\sin, \cos$, in the sense that the derivatives of $a(x,w),b(x,w)$ are certain combination of $a(x,w),b(x,w)$ themselves. All the partials of $y$ is in the form of $c_1aa(-w)+c_2ab(-w)+c_3ba(-w)$, where $c_1,c_2,c_3$ are coefficients consists of $u,w$ and derivative of $u$ (There is no $bb(-w)$ because $bb(-w) = aa(-w)$ by equation (\ref{identitieslast})). Finally if we plugin the all the partials into equation (\ref{equationsofy}), there are terms: $a^3a^3(-w), a^3a^2(-w)b(-w), a^3a(-w)b^2(-w),a^3b^3(-w),b^3a^3(-w)$ etc, the coefficients before every terms will be canceled to 0 using the following relations:
	\begin{gather}
	u_{xx} = 2u^3+xu\\
	u_{xxx} = 6uu_x+u+xu_x
	\end{gather}
    \section*{Acknowledgment}
    I am very grateful to my supervisor Professor Jeremy Quastel for suggesting this problem to me. He gave me many invaluable guidance and discussions on this topic, besides he gave me many important suggestions on my writing of the paper.
	\bibliographystyle{alpha}
	\bibliography{pde_baikrains}
\end{document}